\documentclass[11pt]{amsart}
\usepackage{amscd}
\usepackage{amssymb}
\usepackage[arrow,matrix]{xy}

\newcommand{\C}{{\mathbb C}}
\renewcommand{\L}{{\mathbb L}}
\newcommand{\Q}{{\mathbb Q}}

\newcommand{\V}{{\mathbb V}}
\newcommand{\Vbar}{\overline{\mathbb V}}
\newcommand{\Vtil}{\widetilde{\mathbb V}}
\newcommand{\W}{{\mathbb W}}
\newcommand{\XO}{{X}^\circ}
\newcommand{\Z}{{\mathbb Z}}

\newcommand{\Ga}{{\Gamma}}

\newcommand{\la}{{\lambda}}

\newcommand{\pibar}{\bar{\pi}}
\newcommand{\pitil}{\tilde{\pi}}
\newcommand{\psitil}{\tilde{\psi}}

\newcommand{\A}{{\mathcal A}}
\newcommand{\Abar}{\overline{\mathcal A}}
\newcommand{\CS}{{\mathbb C}^\ast}

\newcommand{\calF}{{\mathcal F}^\bullet}
\newcommand{\calH}{{\mathcal H}}
\newcommand{\calIC}{{\mathcal I \mathcal C}^\bullet}
\newcommand{\M}{{\mathcal M}}
\newcommand{\Mbar}{\overline{\mathcal M}}
\newcommand{\Mtil}{\widetilde{\mathcal M}}
\newcommand{\MO}{{\hspace{-.39cm}{\phantom{\mathcal M}}^{\circ}
		{\hspace{-.1cm}{\mathcal M}}  }}
\newcommand{\calS}{{\mathcal S}^\bul}
\newcommand{\T}{{\mathcal T}}

\newcommand{\tfrak}{{\mathfrak t}}
\newcommand{\teich}{{Teichm\"uller}}

\newcommand{\bul}{\bullet}
\newcommand{\Gr}[2]{{\text{Gr}}^{#1}_{#2}}

\newcommand{\btimes}{\boxtimes}

\newcommand{\Aut}{\operatorname{Aut}}
\newcommand{\Diff}{\operatorname{Diff}}
\newcommand{\Hom}{\operatorname{Hom}}
\newcommand{\Spl}{\operatorname{Sl}}
\newcommand{\Symp}{\operatorname{Sp}}

\newcommand{\ecouple}[6]
{$$
\xymatrix{
#1 \ar[rr]^{#2} &                 & #3 \ar[dl]^{#4} \\
                & #5 \ar[ul]^{#6} &                 }
$$}

\newtheorem{theorem}{Theorem}[section]
\newtheorem{lemma}[theorem]{Lemma}

\newtheorem{corollary}[theorem]{Corollary}
\newtheorem{theoremvoid}{Theorem}

\theoremstyle{definition}

\newtheorem{definition}[theorem]{Definition}

\theoremstyle{remark}

\newtheorem{remark}{Remark}          
	
\newtheorem{notation}{Notation}      
	
\newtheorem{acknowledgements}{Acknowledgements}

\numberwithin{equation}{section}

\begin{document}
\renewcommand{\theequation}{\thesection.\arabic{equation}}
\addtocounter{section}{-1}

\title[The Second Cohomology of the Moduli Space] 
	{The Second Cohomology with Symplectic Coefficients 
	of the Moduli Space of Smooth Projective Curves}

\author{Alexandre I.~Kabanov}

\address{Department of Mathematics\\
Michigan State University\\
Wells Hall \\East Lansing, MI 48824-1027}

\thanks{Research supported in part by an Alfred P.~Sloan Doctoral
Dissertation Fellowship.}

\email{kabanov@math.msu.edu}

\begin{abstract}
Each finite dimensional irreducible rational representation $V$ of the
symplectic group $\Symp_{2g}(\Q)$ determines a generically defined
local system $\V$ over the moduli space $\M_g$ of genus $g$ smooth
projective curves. We study $H^2(\M_g; \V)$ and the mixed Hodge
structure on it. Specifically, we prove that if $g\ge 6$, then the
natural map $IH^2(\Mtil_g; \V) \to H^2(\M_g; \V)$ is an isomorphism
where $\Mtil_g$ is the Satake compactification of $\M_g$. Using the
work of Saito we conclude that the mixed Hodge structure on $H^2(\M_g;
\V)$ is pure of weight $2+r$ if $\V$ underlies a variation of Hodge
structure of weight $r$. We also obtain estimates on the weight of the
mixed Hodge structure on $H^2(\M_g; \V)$ for $3\le g<6$. Results of
this article can be applied in the study of relations in the Torelli
group $T_g$.
\end{abstract}

\maketitle

\section{Introduction}
\label{sec:intro}

The moduli space $\M_g$ of smooth projective curves of genus $g$ is a
quasi-projective variety over $\C$. Its points correspond to
isomorphism classes of smooth projective complex curves of genus
$g$. It has only finite quotient singularities, and therefore behaves
like a smooth variety.

This space has several natural compactifications. In this article we
will be interested in the so called Satake compactification $\Mtil_g$
of $\M_g$. The period map determines an inclusion of $\M_g$ into
$\A_g$, the moduli space of principally polarized abelian
varieties. The Satake compactification $\Mtil_g$ is the closure of
$\M_g$ inside $\Abar_g$, the Satake compactification of $\A_g$. It has
quite complicated singularities at its boundary $\Mtil_g - \M_g$.

Each representation of the algebraic group $\Symp_{2g}$ gives rise to
an orbifold local system over $\M_g$. To explain this we introduce the
mapping class group $\Ga_g$. It is the group of connected components
of the group of the orientation preserving diffeomorphisms of a
compact orientable surface $S$ of genus $g$. The group $\Ga_g$ is the
orbifold fundamental group of $\M_g$, and representations of $\Ga_g$
give rise to orbifold local systems over $\M_g$. There is a natural
surjective map
$$
\Ga_g \rightarrow \Aut (H_1(S; \Z), \cap)
$$
where $\cap$ is determined by the intersection pairing. The right hand
group is isomorphic to $\Symp_{2g}(\Z)$. So each finite dimensional
rational representation $V$ of an algebraic group $\Symp_{2g}$ gives
rise to a {\it symplectic} orbifold local system $\V$ over $\M_g$.

Since $\V$ is generically defined over $\Mtil_g$, one can consider the
intersection cohomology groups $IH^\bul (\Mtil_g; \V)$. There is a
natural restriction map 
$$
IH^\bul (\Mtil_g; \V) \rightarrow H^\bul (\M_g; \V). 
$$
The main result of this article is

\begin{theoremvoid}
(cf. Th.~\ref{thm:main}) The natural restriction map 
$$
IH^k (\Mtil_g; \V) \rightarrow H^k (\M_g; \V)
$$
is an isomorphism when $k=1$ for all $g\ge 3$, and when $k=2$ for all
$g\ge 6$. 
\end{theoremvoid}

The group $H^1(\M_g; \V)$ is easily computed when $g\ge 3$ for all
symplectic local systems $\V$ using Johnson's fundamental work
\cite{johnson} (cf.~\cite{hain:torelli}).

Let $X$ be an algebraic variety. From Saito's work \cite{saito:intro},
\cite{saito:mhm} we know that $H^\bul(X;\V)$ has natural mixed Hodge
structure (MHS) if $\V \to X$ is an admissible polarized variation of
Hodge structure, and $IH^\bul(X;\V)$ has natural mixed Hodge structure
if $\V$ is a generically defined admissible polarized variation of
Hodge structure over $X$. Further if $X$ is compact and $\V$ is pure
of weight $r$, then $IH^\bul(X;\V)$ is pure of weight $k+r$.

\begin{theoremvoid}
(cf. Cor.~\ref{thm:purity}, Cor.~\ref{cor:semipurity}) If
$g\ge 6$ and $\V \rightarrow \M_g$ is a variation of Hodge structure of
weight $r$ whose underlying local system is symplectic, then the
natural mixed Hodge structure on $H^2(\M_g;\V)$ is pure of weight
$2+r$. If $3\le g <6$, then the weights of the mixed Hodge structure
on $H^2(\M_g; \V)$ lie in $\{ 2+r, 3+r \}$. 
\end{theoremvoid}

Each symplectic local system $\V$ associated to an irreducible
representation $V$ of $\Symp_{2g}$ underlies a variation of Hodge
structure over $\M_g$ which is unique up to Tate twist. It is
convenient to fix the weight of the variation of Hodge structure
$\V(\lambda)$ associated to a dominant integral weight $\lambda$. Fix
fundamental weights $\lambda_1, \lambda_2, \dots, \lambda_g$ of
$\Symp_{2g}$. If $\lambda= a_1\lambda_1 + a_2\lambda_2 +\cdots +
a_g\lambda_g$, define $|\lambda|= a_1+ 2a_2+ \cdots +ga_g$. This is
the smallest integer $r$ such that $V(\lambda) \subseteq
H_1(S)^{\otimes r}$. (A good reference is \cite{fulton:harris}.) Then
$\V(\lambda)$ can be realized uniquely as a variation of Hodge
structure of weight $|\lambda|$.

Harer proved in \cite{harer:stability} that the cohomology $H^k(\M_g;
\Z)$ stabilizes when $g \ge 3k$, and Ivanov later improved the range
of stability \cite{ivanov:stability}. \cite{ivanov:twist}. He showed
that $H^k(\M_g; \Z)$ stabilizes when $g\ge 2k+2$. In
\cite{ivanov:twist} Ivanov also proved that $H^k (\M_{g,1};
\V(\lambda))$ is independent of $g$ when $g\ge 2k+2+|\lambda|$. (The
space $\M_{g,1}$ is the moduli space of curves with a marked non-zero
tangent vector.) In \cite{looijenga:stable} Looijenga calculated the
stable cohomology groups of $\M_g$ with symplectic coefficients as a
module over stable cohomology groups of $\M_g$ with trivial
coefficients. In particular, this implies that $H^k (\M_g;
\V(\lambda))$ is independent of $g$ when $g\ge 2k+2+2|\lambda|$.

Looijenga's result also provides very specific information about the
MHS on $H^k (\M_g; \V(\lambda))$. Combined with computations of
$H^k(\M_g; \Q)$ in low dimensions due to Harer \cite{harer:second},
\cite{harer:third}, \cite{harer:fourth}, it implies that $H^k (\M_g;
\V(\lambda))$ is pure of weight $k+|\lambda|$ when $k\le 4$ and $g$ is
in the stability range. In particular, $H^2(\M_g; \V(\lambda))$ is
pure of weight $2+|\lambda|$ when $g\ge 6+2|\lambda|$. Recently,
Pikaart proved in \cite{pikaart} that the stable cohomology $H^k
(\M_g; \Q)$ is pure of weight $k$. Combined with Looijenga's
computations, this shows that $H^k (\M_g; \V(\lambda))$ is pure of
weight $k+|\lambda|$ whenever $g\ge 2k+2+2|\lambda|$. 

Unlike the stability range, our purity range is {\it independent} of
$|\lambda|$. This is important for the following application which was
the motivation for this article.

The Torelli group $T_g$ is the kernel of the surjective homomorphism
$\Ga_g \to \Symp_{2g}(\Z)$. One can consider the Malcev Lie algebra
$\tfrak_g$ associated to $T_g$. (For definitions see
\cite{hain:completion}). This Lie algebra is an analogue of the Lie
algebra associated to the pure braid group on $m$ strings, which is
important in the study of Vassiliev invariants and conformal field
theory. By a result of Johnson \cite{johnson}, $T_g$ is finitely
generated when $g\ge 3$. Thus, $\tfrak_g$ is also finitely generated
when $g\ge 3$. It is not known for any $g\ge 3$ whether $T_g$ is
finitely presented or not.

In \cite{hain:lietor} Hain gives an explicit presentation of
$\tfrak_g$ for $g\ge 3$. More specifically, he proves that for each
choice of $x_0 \in \M_g$ there is a canonical MHS on $\tfrak_g$ which
is compatible with the bracket. Thus,
$$
\tfrak_g \otimes \C \cong \prod_m \Gr{W}{-m} \tfrak_g \otimes \C
$$
where $\Gr{W}{\bul}$ are the graded quotients of the MHS associated to
a choice of $x_0$. Hain proves that for all $g\ge 3$ 
$$
\Gr{W}{\bul} \tfrak_g = \L (H_1 (\tfrak_g)) / (R_g)
$$
where $\L$ stands for the free Lie algebra, and $R_g$ is a set of
relations. According to a result of Johnson \cite{johnson} $H_1
(\tfrak_g)$ is isomorphic as an $\Symp_{2g}$-module to $V(\la_3)$.

Using the above theorem about the MHS on $H^2(\M_g; \V)$ Hain proves
that the relations $R_g$ are quadratic when $g\ge 6$, and quadratic
and possibly cubic when $g=3,4,5$. Moreover, he explicitly calculates
all quadratic relations. This implies that $\tfrak_g$ is finitely
presented for all $g\ge 3$.

We shall outline the proof of the first theorem above. There are three
main steps in the proof. The first step is to notice that if $g\ge 3$,
then the boundary $\Mtil_g - \M_g$ of the Satake compactification has
one irreducible component of codimension two, and all other
irreducible components have codimension three. This immediately
implies that $H^1(\M_g; \V) \cong IH^1(\Mtil_g; \V)$.

The codimension two irreducible component of $\Mtil_g - \M_g$ has a
Zariski open subset isomorphic to $\M_1 \times \M_{g-1}$. We denote it
by $X$. (In the paper we work with a smooth Zariski open subset of
$X$. However this is just a technical detail, and we do not want to
draw an attention to it here.) Let $N^\ast$ be the link bundle of $X$
in $\Mtil_g$. We denote by $\pi$ the corresponding projection. Then
there is an exact sequence
$$
0 \to IH^2(\Mtil_g; \V) \to H^2(\M_g; \V) \to 
H^0(X; R^2 \pi_\ast \V),
$$
and the last morphism factors through the edge homomorphism 
$$
\psi: H^2(N^\ast; \V) \to H^0(X; R^2 \pi_\ast \V)
$$
of the Leray--Serre spectral sequence of $\pi$. Therefore it suffices
to show that $\psi$ is the trivial homomorphism. 

The second step is to understand the link bundle $N^\ast$. Let $L$ be
the pull-back under $pr_2: X \to \M_{g-1}$ of the unit relative
tangent bundle over $\M_{g-1}$, and $\pitil$ be the corresponding
projection $L \to X$. We show that $L$ is a two-to-one unramified
covering of $N^\ast$. (This is done in Section \ref{sec:link}.) Here
we need to assume that $g\ge 4$. Denote by $\Vtil$ the pull-back of
the local system $\V$ to $L$, and by $\psitil$ the edge homomorphism
$H^2(L; \Vtil) \to H^0(X; R^2 \pitil_\ast \Vtil)$ of the Leray--Serre
spectral sequence of $\pitil$. There is a commutative diagram
$$
\begin{CD}
H^2(L; \Vtil) @> \psitil >> H^0(X; 
R^2 \pitil_\ast \Vtil) \\
@AAA  @AAA \\
H^2(N^\ast; \V) @> \psi >>  H^0 (X; 
R^2 \pi_\ast \V) \\
\end{CD}
$$
where both vertical maps are inclusions. This implies that $\psi$ is
trivial, if $\psitil$ is trivial.

The third step is to show that $\psitil$ is trivial. The local system
$\Vtil$ extends to the stratum $X \cong \M_1 \times \M_{g-1}$, and
splits over it according to the branching rule for the standard
inclusion of $\Spl_2 \times \Symp_{2g-2}$ into $\Symp_{2g}$. The
bundle map $\pitil$ respects this splitting. Thus, it suffices to show
that $\psitil$ is trivial for each irreducible symplectic local system
$\Vbar$ over $X$. We complete the computation using Schur's lemma and
the fact, due to Harer \cite{harer:fourth}, that $H^2(\Ga_{g,1};
H^1(S))$ is trivial when $g\ge 4$. (One can also use a result from
\cite[Sec.~7]{harer:third} that $H^2(\Ga_{g,1}; H^1(S))$ is trivial
when $g\ge 9$.)

\begin{acknowledgements}
I would like to thank Richard Hain for his helpful suggestions and
numerous discussions during my graduate studies at Duke University.
John Harer helped to complete the last step in the proof of the main
theorem. I would also like to thank Sloan Foundation for providing me
with the Doctoral Dissertation Fellowship that allowed me to spend the
Fall Semester of 1994 at the Institute for Advanced Study at
Princeton. Conversations with Pierre Deligne, Alan Durfee, Mikhail
Grinberg, Eduard Looijenga, William Pardon, Martin Pikaart, Chad
Schoen, Mark Stern, and Steven Zucker were also very helpful. I would
like to thank the referee who suggested many shortcuts in this paper.
\end{acknowledgements}

\section{Basic facts about the moduli space of curves}
\label{sec:moduli}

In this section we recall the definitions and basic properties of the
moduli spaces of curves, and the corresponding mapping class groups.

The moduli space $\M_{g,r}^s$ parameterizes the isomorphism classes of
smooth complex projective curves of genus $g$ with $s$ marked points
and $r$ marked non-zero holomorphic tangent vectors. The existence of
such moduli spaces follows from geometric invariant theory. These
moduli spaces are known to be normal quasi-projective varieties
\cite[Th.~5.11, Th.~7.13]{mumford:git}.

One can also construct $\M_{g,r}^s$ using \teich\ theory. This approach
allows us to establish the relation between the moduli spaces and the
corresponding mapping class groups.

Let $S$ denote a smooth compact orientable surface of genus $g$. Fix
$s+r$ distinct points $p_1, \dots, p_{r+s}$ on $S$, and $r$ non-zero
tangent vectors $v_1, \dots, v_r$ at points $p_1, \dots, p_r$
respectively. One can consider triples 
$$
(C, (q_1,\dots, q_{r+s}, w_1, \dots, w_r), [f]),
$$
where $C$ is a smooth projective genus $g$ curve,
$q_1, \dots, q_{r+s}$ are distinct points on $C$, $w_1, \dots, w_r$
are non-zero holomorphic tangent vectors at $q_1, \dots, q_r$
respectively, and $f: C \to S$ is an orientation preserving
diffeomorphism such that $f(q_i)= p_i$ and $f_\ast (w_i)= v_i$ (we use
the canonical identification of the holomorphic tangent space with the
underlying real tangent space). We denote by $[f]$ the isotopy class
of $f$ relative to $\{ q_1, \dots, q_{r+s}, w_1, \dots, w_r \}$. Two
triples 
$$
(C_j, (q_1^j, \dots, q_{r+s}^j, w_1^j, \dots, w_r^j), [f_j]),
\quad j=1,2, 
$$
are called {\it equivalent} if there exists a biholomorphism $h: C_1
\to C_2$ such that $h(q_i^1) =q_i^2$, $h_\ast(w_i^1) =w_i^2$, and
$[f_2 \circ h]= [f_1]$ where the homotopy is required to preserve the
marked points and tangent vectors. The space of equivalence classes
$\T_{g,r}^s$ is called the {\it \teich\ space} \cite{harer:review},
\cite[p.~26]{harer:third}. It is known that $\T_{g,r}^s$ is a
contractible complex manifold of dimension $3g-3+s+2r$ when
$2g-2+s+2r>0$.

The {\it mapping class group} $\Ga_{g,r}^s$ is defined to be
$\Diff^+(S)/ \Diff^+_0(S)$, where $\Diff^+(S)$ is the group of
orientation preserving diffeomorphisms of $S$, which leave the marked
points $p_1, \dots, p_{r+s}$ and marked tangent vectors $v_1, \dots,
v_r$ fixed, and $\Diff^+_0(S)$ is the connected component of the
identity. If $g>0$, then the group $\Ga_{g,r}^s$ is torsion free when
either $r>0$, or $s>2g+2$.

The group $\Ga_{g,r}^s$ acts on $\T_{g,r}^s$ as follows. If $g\in
\Ga_{g,r}^s$, then
$$
g (C, (q_1,\dots, w_r), [f])= (C, (q_1,\dots, w_r), [g\circ f]). 
$$
The quotient space $\Ga_{g,r}^s \!\! \setminus \!\! \T^s_{g,r}$ is the
moduli space $\M_{g,r}^s$ of curves with $s$ marked points and $r$
marked tangent vectors. The group $\Ga_{g,r}^s$ acts on $\T_{g,r}^s$
by biholomorphisms, and this action is properly discontinuous and
virtually free. It follows that $\M_{g,r}^s$ is a complex analytic
variety with only finite quotient singularities. This analytic
structure agrees with the one coming from geometric invariant theory.
If $\Ga_{g,r}^s$ is torsion free, then the action is free, and
$\M_{g,r}^s$ is smooth.

\begin{notation}
We shall omit indices $r$ and $s$ from $\T_{g,r}^s$, $\Ga_{g,r}^s$,
and $\M_{g,r}^s$ when they are equal to zero. We shall use both
$\M_1^1$ and $\M_1$ to denote the moduli space of elliptic curves.
\end{notation}

\begin{remark}
One can also consider $\M^{[s]}_{g}$, the moduli space of genus $g$
curves with a marked set of cardinality $s$. It is the quotient of
$\M^s_{g}$ by the natural action of the symmetric group on $s$
letters. This action permutes the marked points.
\end{remark}

The singular locus of $\M_g$ is contained in the locus of curves with
non-trivial automorphisms. When $g\ge 3$, we denote by $\MO_g$ the
locus of curves with only trivial automorphisms. This is a smooth
Zariski open subset of $\M_g$ whose complement has codimension $g-2$.

There are natural surjective morphisms between different moduli spaces
which correspond to forgetting marked points and marked tangent
vectors \cite{knudsen:II}. We will consider the morphisms $\M_g^1 \to
\M_g$ and $\M_{g,1} \to \M_g^1$. The first morphism $\M_g^1 \to \M_g$
is called the ``universal curve''\cite[p.~218]{eisenbud:harris}. Its
fiber over a point $[C]\in \M_g$ is $C/\! \Aut C$. On the level of the
mapping class groups there is a corresponding short exact sequence
\cite{birman:braids}
$$
1 \to \pi_1(S) \to \Ga_g^1 \to \Ga_g \to 1.
$$

The morphism $\M_{g,1} \to \M_g^1$ ``forgets'' the tangent vector, but
remembers its base point. When $g\ge 2$ it is the frame bundle of the
relative holomorphic tangent bundle to the universal curve. On the
level of the mapping class groups there is a corresponding short exact
sequence \cite{birman:braids}
$$
1 \to \Z \to \Ga_{g,1} \to \Ga_g^1 \to 1.
$$  

The composition of the two morphisms discussed above is the morphism
$\M_{g,1} \to \M_g$ obtained by forgetting the tangent vector. If $C$
is a curve without non-trivial automorphisms, then the fiber over $[C]
\in \M_g$ is isomorphic to $T^u C$, the frame bundle of the
holomorphic tangent bundle of the curve $C$. The corresponding
homomorphism of the mapping class groups is $\Ga_{g,1} \to \Ga_g$.

One can also consider finite index level subgroups $\Ga_{g,r}^s [l]$
of $\Ga_{g,r}^s$ for each integer $l$. The level $l$ subgroup is
defined to be the subgroup of $\Ga_{g,r}^s$ which acts trivially on
$H_1(S;
\Z/ l\Z)$. Consequently, one has a short exact sequence
$$
1 \to \Ga_{g,r}^s [l] \to \Ga_{g,r}^s \to
\Symp_{2g}(\Z/ l\Z) \to 1.
$$
The quotient $\Ga_{g,r}^s [l]\setminus \T_{g,r}^s$ is isomorphic to
$\M_{g,r}^s [l]$, the moduli space of smooth projective curves with a
level $l$ structure which is defined in Section \ref{sec:comp}. 

It is well-known that for all $g\ge 1$ and $l\ge 3$, the group
$\Ga_{g,r}^s [l]$ acts freely on $\T_{g,r}^s$. Thus for each $l\ge 3$
the moduli space $\M_{g,r}^s [l]$ is a smooth finite cover of
$\M_{g,r}^s$.

When $\M_{g,r}^s$ is different from $\M_1$ and $\M_2$ each
representation of $\Ga_{g,r}^s$ determines an orbifold local system
over $\M_{g,r}^s$. When $\M_g$ is either $\M_1$ or $\M_2$ we consider
only such representations of $\Ga_g$ that for each $[C] \in \M_g$
represented by a curve with only two automorphisms, the stabilizer of
$(C,[f])\in \T_g$ acts trivially on the representation space. These
representations give rise to orbifold local systems over $\M_1$ and
$\M_2$.

Let $V$ be a representation of $\Ga_{g.r}^s$ on a rational vector
space, and let $\V$ be the associated orbifold local system over
$\M_{g,r}^s$. The contractibility of the \teich\ space implies that
for all $g\ge 1$:

$$
H^\bul (\Ga_{g,r}^s; V) \cong H^\bul (\M_{g,r}^s; \V) \cong
H^\bul (\M_{g,r}^s [l]; \V[l])^{\Symp_{2g} (\Z/ l\Z)}.
$$


\section{Compactifications of the moduli space of curves}
\label{sec:comp}

In this section we recall some basic properties of the Satake
compactification and the Deligne--Mumford compactification of the
moduli spaces of curves.

We start with the Deligne--Mumford compactification of $\M_g^s$. A
{\it stable curve} is a reduced connected curve which has only nodes
as singularities, and a finite automorphism group
\cite{deligne:mumford}. The {\it Deligne--Mumford compactification}
$\Mbar_g^s$ of $\M_g^s$ is the moduli space of stable projective
curves. It is a normal projective variety in which $\M_g^s$ is a
Zariski open subset \cite{deligne:mumford},
\cite[Th.~5.1]{mumford:stable}. The singularities of $\Mbar_g^s$ are
contained in the locus of stable curves with non-trivial automorphisms
\cite[p.~218]{eisenbud:harris}.

We will describe the boundary $\Mbar_g^s - \M_g^s$ in the case when
$s=0$. The boundary $\Mbar_g-\M_g$ is the union of irreducible
divisors
$$
\bigcup_{i=0}^{[{g/ 2}]} \Delta_i,
$$
where each divisor $\Delta_i$ has the following property. When $i=0$
there is birational morphism $\Mbar_{g-1}^{[2]} \to \Delta_0$; when
$1\le i < g-i$ there is birational morphism $\Mbar_i^1 \times
\Mbar_{g-i}^1 \to \Delta_i$; and when $i=g-i$ there is a birational
morphism from the $\Z/2\Z$-quotient of $\Mbar_i^1 \times \Mbar_i^1$ to
$\Delta_i$.

\begin{definition}
(cf. \cite[Def.~10.5]{popp}) A {\it level $l$ structure} on a
stable curve $C$ is a symplectic monomorphism $H^1(C; \Z/l\Z) \to
(\Z/l\Z)^{2g}$, where $(\Z/l\Z)^{2g}$ has the standard symplectic
structure.
\end{definition}

Note that a level $l$ structure on a smooth curve $C$ is just a choice
of a symplectic basis for $H^1(C;\Z/l\Z)$, or, equivalently, for
$H_1(C;\Z/l\Z)$ because the symplectic form determines the canonical
identification between homology and cohomology. The same is true for a
singular stable curve $C$ whose dual graph is a tree. 

From now on we assume that $l\ge 3$. Denote by $\M_g [l]$ the moduli
space of smooth curves with a level $l$ structure. It is isomorphic to
the quotient of $\T_g$ by the action of $\Ga_g [l]$ (cf.
Sec.~\ref{sec:moduli}). The moduli space $\M_g [l]$ is a smooth
quasi-projective variety, and the forgetful morphism $\M_g[l] \to
\M_g$ is a Galois covering \cite[Prop.~5.8]{deligne:mumford},
\cite[Th.~1.8]{oort:steenbrink}.

When $g\ge 2$ and $l\ge 3$ there exists the moduli space of stable
curves with a level $l$ structure $\Mbar_g [l]$, which is a
compactification of $\M_g [l]$ \cite[p.~106]{deligne:mumford},
\cite[Bem.~1]{mostafa}, \cite[Rem.~2.3.7]{jong:pikaart}. This is a
projective variety according to \cite[Th.~4, III.8]{mumford:red}, and
there is a finite morphism $\Mbar_g [l] \to \Mbar_g$ determined by
forgetting a level $l$ structure.

In \cite{mostafa} Mostafa proves that $\Mbar_g [l]$ is not smooth, at
least when $g\ge 3$. However, in this article we are interested in
particular strata of the boundary of $\Mbar_g [l]$. The irreducible
component $\Delta_1$ of the boundary of $\Mbar_g$ contains a Zariski
open subset isomorphic to $\M_1^1 \times \M_{g-1}^1$. Consider the
inverse image of this subset under the finite morphism above. It is a
finite disjoint union of locally closed subvarieties of codimension
one each of which is isomorphic to $\M_1^1 [l] \times \M_{g-1}^1 [l]$.
According to \cite[Lem.~1]{mostafa}, \cite[p.~240]{looijenga:prym} all
points of this inverse image are smooth points of $\Mbar_g [l]$.

To introduce the Satake compactification $\Mtil_g$ of $\M_g$ we use
the space $\A_g$, the moduli space of principally polarized abelian
varieties of dimension $g$. It is the quotient of the Siegel
upper-half space by the action of $\Symp_{2g}(\Z)$. The space $\A_g$
is a quasi-projective variety \cite[Th.~7.10]{mumford:git}. Among
other compactifications, it admits the Satake compactification
$\Abar_g$ which is a projective variety \cite{satake}.

The moduli space $\M_g$ is isomorphic to the image of the period map
$\M_g \to \A_g$ which is a locally closed subvariety of $\A_g$
\cite[Cor.~3.2]{oort:steenbrink}. The closure $\Mtil_g$ of $\M_g$ in
the Satake compactification $\Abar_g$ of $\A_g$ is called the {\it
Satake compactification of} $\M_g$ (cf. \cite{baily:1962}). There
exists a birational morphism $\alpha: \Mbar_g \to \Mtil_g$ which is
the identity on $\M_g$, and sends the point $[C]$ corresponding to a
stable curve $C$ to the polarized Jacobian of its normalization
\cite[p.~211]{knudsen:III}.

The image of the boundary $\Mbar_g - \M_g$ under $\alpha$ is the
boundary $\Mtil_g - \M_g$ of the Satake compactification. It follows
that when $g\ge 3$ the boundary $\Mtil_g - \M_g$ has $[g/2]$
irreducible components each of which except one has codimension three
in $\Mtil_g$. The irreducible component $\Phi_1$ which is the image of
$\Delta_1 \subset \Mbar_g$ has codimension two. It contains a Zariski
open subset isomorphic to $\M_1 \times \M_{g-1}$.

One can also construct the {\it Satake compactification} $\Mtil_g [l]$
of $\M_g [l]$. Denote by $\A_g [l]$ the moduli space of principally
polarized abelian varieties with a level $l$ structure. A point in
$\A_g [l]$ is represented by an abelian variety $A$ of dimension $g$
and a symplectic basis of $H_1 (A; \Z/l\Z)$. It is a quasi-projective
variety \cite[Th.~7.9]{mumford:git}, \cite[Th.~1.8]{oort:steenbrink}
which is smooth when $l\ge 3$. The space $\A_g [l]$ has the Satake
compactification $\Abar_g [l]$ which is a normal projective variety
\cite[p.~124]{popp}, \cite{satake}.

If $g=1,2$, then $\M_g [l]$ is isomorphic to a Zariski open subset of
$\A_g [l]$, and we define the Satake compactification $\Abar_g [l]$ of
$\A_g [l]$ to be the Satake compactification of $\M_g [l]$. If $g\ge 3$,
then the morphism $\M_g [l] \to \A_g [l]$ is not injective. In this
case we define $\Mtil_g [l]$ to be the normalization of $\Mtil_g$ with
respect to $\M_g [l]$.

It follows from this definition that $\Mtil_g [l]$ is a projective
variety \cite[III.8, Th.~4]{mumford:red}, and that the morphism
$\M_g [l] \to \M_g$ extends to a finite morphism $\Mtil_g [l]
\to \Mtil_g$. One can also show that there is a birational morphism
$\alpha^l:\Mbar_g [l] \to \Mtil_g [l]$ with connected fibers which is
the identity on $\M_g [l]$, and fits into the commutative diagram
$$
\begin{CD}
\Mbar_g [l] @>\alpha^l >> \Mtil_g [l] \\
@VVV                       @VVV \\
\Mbar_g     @>\alpha >>  \Mtil_g. \\
\end{CD}
$$

The boundary $\Mtil_g [l] - \M_g [l]$ is the union of irreducible
components each of which has codimension either two, or three in
$\Mtil_g [l]$. The image of each component $\Phi_1^\beta$ of
codimension two under the morphism $\Mtil_g [l] \to \Mtil_g$ is the
codimension two component $\Phi_1$ of $\Mtil_g - \M_g$. One can show
that each $\Phi_1^\beta$ contains a Zariski open subset $Z_\beta$ such
that these subsets do not intersect each other, and each of them is
isomorphic to a smooth Zariski open subset of $\M_1 [l] \times
\M_{g-1} [l]$.


\section{Codimension two stratum of the Satake compactification}
\label{sec:link}

In this section we analyze the link of the codimension two boundary
stratum $\Phi_1$ inside the Satake compactification of the moduli
space $\Mtil_g$. More precisely, we study the local links of the
points in a smooth Zariski open subset of $\Phi_1$, and we show that
$\Mtil_g$ is equi-singular along this Zariski open subset. We will
need this in Section \ref{sec:main}. For the rest of this section we
assume that $g\ge 4$.

Recall that $\Phi_1$ contains a Zariski open subset $X$ isomorphic to
$\M_1 \times \M_{g-1}$. We identify it with $\M_1 \times \M_{g-1}$.
Then a point in $X$ is represented by a pair of isomorphism classes of
curves $([C_1],[C_2])$. Let $\XO$ be a Zariski open subset of $X$
defined as follows. Recall that in Section \ref{sec:moduli} we defined
$\MO_g$ to be the locus of curves with only trivial automorphisms when
$g\ge 3$. We define $\MO_1$ to be the locus of elliptic curves with
exactly two automorphisms. Then $\XO$ is the subset of $X$
corresponding to $\MO_1 \times \MO_{g-1}$. In this section we study
the link of $\XO$ in $\M_g \cup \XO \subset \Mtil_g$.

Let $N$ be a regular neighborhood of $\XO$ in $\M_g \cup \XO$. The
complement $N^\ast= N- \XO$ is a deleted regular neighborhood of
$\XO$.

Recall that $\M_{g-1,1} \to \M_{g-1}$ is a surjective morphism defined
by forgetting the holomorphic tangent vector. Let $L_2$ be the inverse
image of $\MO_{g-1}$ in $\M_{g-1,1}$, and $\pitil_2: L_2 \to
\MO_{g-1}$ be the corresponding map. The fiber of $\pitil_2$ over
$[C_2] \in \MO_{g-1}$ is $T^u C_2$, the punctured holomorphic tangent
bundle. Denote by $L$ the product $\MO_1 \times L_2$, and by $\pitil$
the pull-back of $\pitil_2$ to $\XO$:
$$
\begin{CD}
L= \MO_1 \times L_2 @> pr_2 >> L_2 \\
@V \pitil VV             @V \pitil_2 VV \\
\XO= \MO_1 \times \MO_{g-1} @> pr_2 >> \MO_{g-1}. \\ 
\end{CD}
$$

\begin{lemma}
\label{lem:link}
The bundle $\pitil: L= \MO_1 \times L_2 \to \XO$ is a two-to-one
unramified cover of the punctured regular neighborhood $N^\ast$. The
corresponding fix point free action of $\Z/2\Z$ on $L$ sends a vector
$v$ to $-v$.
\end{lemma}

\begin{proof}
The morphism $\M_{g-1,1} \to \M_{g-1}$ factors as 
$$
\M_{g-1,1} \to \M_{g-1}^1 \to \M_{g-1}.
$$
Denote by $Y_2$ the inverse image of $\MO_{g-1}$ under the second
morphism. Then the commutative diagram above factors as
$$
\begin{CD}
L= \MO_1 \times L_2 @> pr_2 >> L_2 \\
@V \pi^c VV             @V \pi^c_2 VV \\
Y= \MO_1 \times Y_2 @> pr_2 >> Y_2 \\
@V \pibar VV              @V \pibar_2 VV \\
\XO= \MO_1 \times \MO_{g-1} @> pr_2 >> \MO_{g-1}, \\ 
\end{CD}
$$
where $\pi_2^c$ (resp. $\pibar_2$) is the restriction of $\M_{g-1,1}
\to \M_{g-1}^1$ (resp. $\M_{g-1}^1 \to \M_{g-1}$) to $L_2$
(resp. $Y_2$), and $\pi^c$ (resp. $\pibar$) is its pull-back along
$pr_2$.

At the same time $Y= \MO_1 \times Y_2$ is isomorphic to a smooth
Zariski open subset of the boundary component $\Delta_1$ in the
Deligne--Mumford compactification. We identify $Y$ with this Zariski
open subset. Then the morphism $\pibar: Y \to \XO$ is the restriction
of the morphism $\alpha: \Mbar_g \to \Mtil_g$ to $Y$.

The morphism $\alpha$ is the identity when restricted to $\M_g$.
Therefore a deleted regular neighborhood $N^\ast$ of $\XO$ in $\M_g
\cup \XO$ and a deleted regular neighborhood of the divisor $Y$
in $\M_g \cup Y \subset \Mbar_g$ can be chosen to be the same.

The deleted neighborhood of $Y$ is homeomorphic to the punctured
normal bundle of $Y$ in $\M_g \cup Y$. Note that the
only non-trivial automorphism of a pair $(C_1,x_1),(C_2,x_2)$
representing a point in $Y$ is induced by the elliptic involution of
$(C_1,x_1)$. It follows
that $\Z/2\Z$ acts on the space of versal deformations of the stable
curve $(C_1,x_1),(C_2,x_2)$, and this action fixes the divisor that is
the locus of the singular curves
\cite[Chap.~13, Lem.~(1.6)]{acgh:II}. Therefore the fiber of the
normal bundle of $Y \subset \Delta_1$ at the point
$[(C_1,x_1),(C_2,x_2)]$ is isomorphic to the $\Z/2\Z$ quotient of
$T_{x_1}C_1 \otimes T_{x_2}C_2$, where the generator of $\Z/2$ acts as
$-id$. Thus $N^\ast$ is the $\Z/2\Z$ quotient of the $\CS$-bundle $L'$
over $Y$ whose fiber at $[(C_1,x_1),(C_2,x_2)]$ is $T_{x_1}C_1 \otimes
T_{x_2}C_2 - \{ 0 \}$.

It is well-known that the moduli space of elliptic curves $\M_1$ is
isomorphic to $\C$. It contains two distinguished points that
correspond to the two elliptic curves with exceptional
automorphisms. It follows that the space $\MO_1$ is isomorphic to
$\C- \{ 2\ \text{points} \}$. All line bundles over this space are
trivial. Therefore the bundle $L'$ is the pull-back of the punctured
relative tangent bundle of the morphism $Y_2 \to \MO_{g-1}$.

The punctured relative tangent bundle of the morphism $Y_2 \to
\MO_{g-1}$ is $\pi^c_2: L_2 \to Y_2$. Hence, one has a commutative
diagram
$$
\begin{CD}
L' @>>> L_2 \\
@VVV  @V \pi^c_2 VV \\
Y=\MO_1 \times Y_2 @> pr_2 >> Y_2, \\
\end{CD}
$$
where $L'$ is the pull-back of $L_2$. We conclude that the bundles
$L'$ and $L$ are isomorphic, and $N^\ast$ is the $\Z/2\Z$ quotient of
$L$, where $\Z/2\Z$ action sends a vector in a fiber of $\pi^c$ to its
opposite.
\end{proof}

It follows from the lemma above that $\Mtil_g$ is equi-singular along
$\XO$. We expressed $N^\ast$ as a bundle $\XO$ whose fiber over the
point $([C_1,C_2])$ is equal to $T^u C_2$, the frame bundle of the
holomorphic tangent bundle of $C_2$. 


\section{Main theorem}
\label{sec:main}

In this section we prove the main theorem of this article. The proof
consists of a sequence of lemmas and propositions. We assume that the
reader is familiar with intersection cohomology, and suggest the
references \cite{bbd}, \cite{borel:ic}, \cite{gmII}.

\begin{notation}
For the rest of the paper we omit $R^\bul$ from the notation for the
derived functors. For example, if $f:X\to Y$ is a continuous map
between topological spaces, then $f_\ast = R^\bul f_\ast$.
\end{notation}

As we mentioned before each representation of the mapping class group
$\Ga_g$, at least when $g\ge 3$, determines an orbifold local system
over $\M_g$. In this section we consider only the {\it symplectic}
local systems, that is local systems arising from finite dimensional
rational representations of the algebraic group $\Symp_{2g}$. We fix a
symplectic representation $V$ of $\Ga_g$, and denote the corresponding
orbifold local system by $\V$.

\begin{theorem}
\label{thm:main}
The natural map $IH^k(\Mtil_g;\V) \to H^k(\M_g;\V)$ induced by the
inclusion is an isomorphism, when
\begin{center}
$\begin{array}{l@{,\quad}l}
k=0 & g\ge 1;\\
k=1 & g\ge 3; \\
k=2 & g\ge 6.
\end{array}$
\end{center}
\end{theorem}

The first statement is trivial and included only for the sake of
completeness. The statement concerning the first cohomology is also
rather simple. Indeed, in Section \ref{sec:comp} we saw that if $g\ge
3$, then the boundary $\Mtil_g - \M_g$ of the Satake compactification
has codimension two in $\Mtil_g$. This, and the properties of
intersection cohomology immediately imply the statement of the theorem
for $k=1$. The non-trivial part of this theorem concerns the second
cohomology.

\begin{remark}
If $g\ge 3$, then the map $IH^1(\Mtil_g;\V) \to H^1(\M_g;\V)$ is an
isomorphism for an arbitrary orbifold local system $\V$ determined by
a representation of $\Ga_g$ on a rational vector space. This can be
easily seen from the above argument.
\end{remark}

Combining this with the computations of $H^1(\M_g;\V)$ in
\cite{hain:torelli}, \cite{johnson} one gets the following corollary.

\begin{corollary}
If $g\ge 3$ and $\V(\lambda)$ is a generically defined local system
corresponding to the representation of $\Symp_{2g}$ with the highest
weight $\lambda$, then
$$
IH^1 (\Mtil_g; \V(\lambda)) \cong 
\left \{ \begin{array}{ll}
\Q & \quad {\mbox{when}} \quad \lambda= \lambda_3 ; \\
0  & \quad {\mbox{otherwise}}. \qquad \qquad \qed
\end{array} \right.
$$ 
\end{corollary}

The rest of this section is devoted to the proof of the isomorphism in
second cohomology. We assume that $g\ge 4$. Recall that we denote by
$\Phi_1$ the codimension two irreducible component of the boundary of
$\Mtil_g$, and by $\XO$ its Zariski open subset isomorphic to $\MO_1
\times \MO_{g-1}$.

\begin{notation}
We denote by $\calS$ the intersection cohomology sheaf $\calIC (\V)$
on $\Mtil_g$ corresponding to the local system $\V$. The following
diagram defines the notation for the inclusions:
$$
\M_g \stackrel{i}{\hookrightarrow} \M_g \cup \XO 
\stackrel{j}{\hookleftarrow} \XO.
$$
\end{notation}

First, we use again that the boundary of $\Mtil_g$ has only one
irreducible component of codimension two, namely $\Phi_1$, and all
other irreducible components have codimension three. This and the
properties of intersection cohomology imply that the restriction
$$
IH^2(\Mtil_g; \V) \to IH^2(\M_g \cup \XO; \V)
$$ 
is an isomorphism, and there is an exact sequence
\begin{multline*}
0 \to IH^2(\M_g \cup \XO; \V) \to H^2(\M_g; \V) \stackrel{\phi}{\to}
\\ H^3 (\XO; j^! \calS) \cong H^0 (\XO; \calH^3 j^! \calS).
\end{multline*}
Therefore to prove the theorem it suffices to show that $\phi$ from
the exact sequence above is the zero morphism.

The distinguished triangle 
\ecouple{j^! \calS} {}
	{j^* \calS} {}
	{j^* i_* \V} {[1]} 
implies that $\calH^3 j^! \calS \cong \calH^2 j^* i_* \V$. Then the
morphism $\phi$ composed with this isomorphism can be factored as
$$
H^2(\M_g; \V) \to H^2(\XO; j^* i_* \V) \stackrel{\psi}{\to} 
H^0(\XO; \calH^2 j^* i_* \V).
$$
The sheaf $j^* i_* \V$ is called the {\it local link cohomology
functor} \cite[p.~57]{durfee:saito}. It expresses the cohomology of
$N^\ast$, the link of $\XO$ in $\M_g \cup \XO$. We denote by $\pi$ the
corresponding projection $N^\ast \to \XO$. Then the morphism $\psi$
from the exact sequence above can be written as
$$
\psi: H^2(N^\ast; \V) \to H^0(\XO; \calH^2 \pi_\ast \V). 
$$
One can easily check that $\psi$ is the edge homomorphism associated
to the Leray--Serre spectral sequence determined by $\pi$. 

In order to prove the theorem it is enough to show that $\psi$ is the
trivial homomorphism when $g\ge 6$, and the rest of this section deals
with the proof of this fact.

First we want to understand the behavior of the local system $\V$ over
$N^\ast$. We start with the following lemma. 

\begin{lemma}
\label{lem:split}
The orbifold local system $\V$ over $N^\ast$ splits into a direct sum
of symplectic orbifold local systems determined by rational
representations of $\Spl_2 \times \Symp_{2g-2}$. 
\end{lemma}

\begin{proof}
Recall that a symplectic orbifold local system is determined by a
representation of $\Ga_g$ which is the pull-back of an algebraic
representation $V$ of $\Symp_{2g}$. 

Choose a level $l\ge 3$. The inverse image of $\XO \subset \Mtil_g$ in
$\Mtil_g [l]$ has several connected components. Let $N_l^\ast$ be a
deleted regular neighborhood of one of them. Then one has a
commutative diagram
$$
\begin{CD}
N_l^\ast @>>> \M_g [l] @>>> \A_g [l] \\
@VVV            @VVV        @VVV  \\
N^\ast @>>>    \M_g     @>>> \A_g.  \\
\end{CD}
$$
(Recall that $A_g$ stands for the moduli space of principally
polarized abelian varieties.) Denote by $\V_l$ the pull-back of $\V$
to $\M_g [l]$, and by $\V_l'$ the local system over $A_g [l]$
determined by $V$. Both $\V_l$ and $\V_l'$ are genuine local systems,
and $\V_l$ is the pull-back of $\V_l'$ under $\M_g [l] \to \A_g [l]$.

The product $\A_1 \times \A_{g-1}$ is canonically embedded in $\A_g$.
Its inverse image under $\A_g [l] \to \A_g$ consists of several
connected component, each of which is isomorphic to $\A_1 [l]\times
\A_{g-1} [l]$. The image of $N_l^\ast$ in $\A_g [l]$ is contained in a
tubular neighborhood of one of these connected components. The local
system $\V_l'$, restricted to this connected component, splits
according to the branching law of the inclusion $\Spl_2 \times
\Symp_{g-2} \hookrightarrow \Symp_{2g}$. It follows that the local
system $\V_l$ splits over $N_l^\ast$ according to the same branching
law. In addition, $\V_l$ is constant on the fibers of the composite
$$
N_l^\ast \longrightarrow N^\ast \stackrel{\pi}{\longrightarrow}
\XO.
$$
Thus the splitting of $\V_l$ over $N_l^\ast$ descends to the splitting
of $\V$ over $N^\ast$.
\end{proof}

Our aim is to show that the morphism $\psi$ is trivial. Therefore
without loss of generality we can assume that $\V$ is a local system
over $N^\ast$ determined by an irreducible algebraic representation of
$\Spl_2 \times \Symp_{2g-2}$ with highest weight $(\mu, \nu)$. Note
that $\mu$ is just a non-negative integer.

We consider two cases. First, assume that $\mu$ is odd. The morphism
$N_l^\ast \to N^\ast$ is a Galois covering with the Galois group
$\Spl_2 (\Z/l\Z) \times \Symp_{2g-2} (\Z/l\Z)$. The element $(-id,
id)$ of this group leaves the fibers of $N_l^\ast \to N^\ast \to
\XO$ fixed because it corresponds to the involution of the
elliptic curve, and acts as $-id$ on the local system $\V_l$. It
follows that $\calH^2 \pi_\ast \V$ is the trivial local system, and we
have nothing more to prove.

Next, assume that $\mu$ is even. Then $(-id, id)$ acts trivially on
the local system $\V_l$, and therefore the local system $\V$ extends
to $\XO$. This means that $\V$ is the restriction to $N^\ast$ of a
local system defined on the whole regular neighborhood $N$ of $\XO$.
Denote the restriction of this local system to $\XO$ by $\Vbar$. Then
$\Vbar$ is isomorphic to $\W_1 (\mu) \btimes \W_2 (\nu)$, the
symplectic local system over $\MO_1 \times \MO_{g-1}$ determined by
the highest weight $(\mu, \nu)$. The local system $\V$ is the
pull-back of $\Vbar$ under $\pi: N^\ast \to \XO$.

Lemma \ref{lem:link} says that $N^\ast$ is the $\Z/2\Z$ quotient of
the bundle $L$ defined in Section \ref{sec:link}, and $\pi$ is induced
by the projection $\pitil: L \to \XO$. We denote by $\Vtil$ be the
pull-back of $\V$ to $L$. Then
$$
\Vtil =\pitil^\ast \Vbar \cong\pitil^\ast (\W_1(\mu) \btimes \W_2(\nu)).
$$

Let $B_\bul$ be the Leray--Serre spectral sequence determined by
$\pi$, and $A_\bul$ be the Leray--Serre spectral sequence determined
by $\pitil$. Let $\psitil$ be the edge homomorphism $H^2(L; \Vtil) \to
H^0 (\XO;\calH^2 \pitil_\ast\Vtil)$ associated $A_\bul$. The two-fold
covering map $L \to N^\ast$ induces the map of the spectral sequences
$B_\bul \to A_\bul$. One has for each $q$ \cite[p.~85]{bredon:sheaf}
$$
\calH^q \pi_\ast \V = (\calH^q \pitil_\ast \Vtil)^{\Z/2\Z}. 
$$
It follows that the induced map $B_2^{0,q} \to A_2^{0,q}$ is an
inclusion of global $\Z/2\Z$ invariants. The homomorphism $H^2(N^\ast;
\V) \to H^2(L; \Vtil)$ is also an inclusion of $\Z/2\Z$ invariants,
and one has a commutative diagram
$$
\begin{CD}
H^2(L; \Vtil) @>\psitil >> H^0 (\XO;\calH^2 \pitil_\ast\Vtil)\\
@AAA                           @AAA                       \\
H^2(N^\ast; \V) @>\psi >> H^0(\XO;\calH^2 \pi_\ast\V), \\
\end{CD}
$$
where both vertical maps are inclusions. It follows that if $\psitil$
is trivial, then $\psi$ is trivial. We shall show that $\psitil$ is
trivial.

Note that $\pitil_\ast \Vtil$ is quasi-isomorphic to $\pitil_\ast \Q
\otimes \Vbar$. It follows that $ \calH^2 \pitil_\ast \Vtil$ is
isomorphic to $\calH^2 \pitil_\ast \Q \otimes \Vbar$.

\begin{lemma}
The local system $\calH^2 \pitil_\ast\Q$ over $\XO$ is
isomorphic to $\W_1(0) \btimes \W_2(\nu_1)$.
\end{lemma}

\begin{proof}
The bundle $\pitil$ is the pull-back of the bundle $\pitil_2: L_2 \to 
\MO_{g-1}$ to $\XO$ (see Lem.~\ref{lem:link}). It follows that
the local system $\calH^2 \pitil_\ast \Q$ is the exterior tensor
product of the constant local system $\W_1(0)$ over $\MO_1$ and
$\calH^2 \pitil_{2\ast} \Q$.

Recall that $\pitil_2$ factors as 
$$
L_2 \stackrel{\pi^c_2}{\longrightarrow} Y_2 
\stackrel{\pibar_2}{\longrightarrow} \MO_{g-1},
$$
where $\pibar_2$ is the restriction of the universal curve to $Y_2$,
and $\pi_2^c$ is a punctured relative tangent bundle to
$\pibar_2$. Therefore we have a Gysin long exact sequence of local
systems
\begin{equation}
\label{gysin}
\cdots \to  \calH^0 \pibar_{2\ast} \Q \stackrel{e}{\longrightarrow} 
\calH^2 \pibar_{2\ast} \Q \to  \calH^2 \pitil_{2\ast} \Q 
\to  \calH^1 \pibar_{2\ast} \Q \to  0
\end{equation}
where $e$ is the multiplication by the Euler class. The Euler class is
non-zero, because the genus of $C_2$ is greater than one. It follows
that $e$ is an isomorphism on rational cohomology. Thus we conclude
that $\calH^2 \pitil_{2\ast} \Q$ is isomorphic to $\calH^1
\pibar_{2\ast} \Q$.

The local system $\calH^1 \pibar_{2\ast} \Q$ is isomorphic to
$\W_2(\nu_1)$, the local system corresponding to the standard
representation of $\Symp_{2g-2}$. It follows that
$$
\calH^2 \pi^\ast\Q \cong \W_1(0)\btimes \W_2(\nu_1). \quad \qed
$$
\renewcommand{\qed}{}
\end{proof}

The lemma above shows that the edge homomorphism $\psitil$ is of the
form
\begin{equation*}
H^2(L;\pitil^\ast(\W_1(\mu)\btimes\W_2(\nu))) \to  
H^0(\XO;\W_1(\mu)\btimes (\W_2(\nu_1)\otimes\W_2(\nu))).
\end{equation*}

\begin{lemma}
The space $H^0(\XO;\W_1(\mu)\btimes (\W_2(\nu_1)\otimes\W_2(\nu)))$ is
isomorphic to $\Q$ if $\mu=0$ and $\nu=\nu_1$, and zero otherwise.
\end{lemma}

\begin{proof}
Applying the K\"unneth formula one gets
\begin{align*}
H^0(\XO;\W_1(\mu)\btimes
(\W_2(\nu_1)&\otimes\W_2(\nu))) \cong \\ H^0(\MO_1;\W_1(\mu))&\otimes
H^0(\MO_{g-1};\W_2(\nu_1)\otimes\W_2(\nu)).
\end{align*}
The zero cohomology of a space with coefficients in a local system is
equal to the space of global invariants of the local system. An
irreducible symplectic local system has no global invariants unless it
is constant. This implies that $H^0(\MO_1; \W_1(0)) \cong \Q$, and
$H^0(\MO_1; \W_1(\mu))=0$ if $\mu\ne 0$.

Similarly, $H^0(\MO_{g-1};\W_2(\nu_1)\otimes\W_2(\nu))$ is equal to
zero, unless the local system $\W_2(\nu_1)\otimes\W_2(\nu)$ contains a
constant local system as a direct summand. This occurs if and only if
the tensor product $W_2(\nu_1)\otimes W_2(\nu)$ of irreducible
representations of $\Symp_{2g-2}(\Q)$ contains a copy of the trivial
representation. It is known that all irreducible representations of
the symplectic group are self-dual. Therefore the trivial part of that
representation is equal to
\begin{align*}
(W_2(\nu_1)\otimes W_2(\nu))^{\Symp_{2g-2}(\Q)}  & =
W_2(\nu_1)\otimes_{\Symp_{2g-2}(\Q)}W_2(\nu)   \\ & \cong
\Hom_{\Symp_{2g-2}(\Q)}(\W_2(\nu_1);\W_2(\nu)).
\end{align*}
By Schur's lemma the latter term is isomorphic to $\Q$, if
$\nu=\nu_1$, and $0$ otherwise.
\end{proof}

It follows that $\psitil$ is trivial unless $\Vbar \cong \W_1(0) \btimes
\W_2(\nu_1)$.  In the remaining part of this section we study this
case. To simplify the notation we denote $\W_1(0) \btimes \W_2(\nu_1)$
by $\W_2(\nu_1)$.

\begin{lemma}
\label{lem:g6}
If $g\ge 6$, then the homomorphism
$$ 
\psitil: H^2(L;\pitil^\ast\W_2(\nu_1)) \to  
H^0(\XO; \W_2(\nu_1)^{\otimes 2})
$$ 
is the zero map.
\end{lemma}

\begin{proof}
Note that $\psitil$ factors through the $A_\infty^{0,2}$ term of the
spectral sequence
$$
A_2^{p,q} = H^p(\XO; 
\calH^q \pitil_\ast \Q \otimes \W_2(\nu_1)) \Rightarrow
H^{p+q} (L; \pitil^\ast \W_2(\nu_1)). 
$$
Therefore, it suffices to prove that
$A_\infty^{0,2}=0$. 

The morphism $\pitil_2: L_2 \to  \MO_{g-1}$ gives rise to a
Leray--Serre spectral sequence
\begin{equation*}
C_2^{p,q}=H^p(\MO_{g-1};\calH^q\pitil_{2\ast} \Q 
\otimes \W_2(\nu_1)) \Rightarrow
H^{p+q}(L_2;\pitil_{2}^\ast \W_2(\nu_1)).
\end{equation*}
The bundle $\pitil: L \to \XO$ is the pull-back of $L_2$, and the
local system $\W_2(\nu_1)$ over $\XO$ is also the pull-back from the
second factor. It follows that the morphism of spectral sequences
$C_\bul \to A_\bul$ induced by the projection $pr_2: \XO \to
\MO_{g-1}$ is an inclusion of a direct summand. (Here we mean that
for each $(r,p,q)$ the term $C_r^{p,q}$ is a direct summand of
$A_r^{p,q}$, and all differentials $d_r$ respect this splitting.)

Note that $A_2^{0,2} \cong C_2^{0,2}$. Indeed,
\begin{align*}
A_2^{0,2}=H^0(\XO; \W_2(\nu_1)^{\otimes 2}) &\cong 
H^0(\MO_1;\Q)\otimes H^0(\MO_{g-1}; \W_2(\nu_1)^{\otimes 2})\\
&\cong H^0(\MO_{g-1};\W_2(\nu_1)^{\otimes 2})=C_2^{0,2}
\end{align*}
because $\calH^2 \pitil_{2\ast} \Q \cong \W_2(\nu_1)$ according to exact
sequence \eqref{gysin}. It follows that $A_{\infty}^{0,2} \cong
C_{\infty}^{0,2}$.

The final step is to show that $C_\infty^{0,2}=0$. There is a
surjective homomorphism
$$
H^2(L_2; \pitil_2^\ast\W_2(\nu_1)) \to C_\infty^{0,2},
$$
associated to the spectral sequence $C_\bul$. Therefore it suffices
to show that 
$$
H^2(L_2; \pitil_2^\ast \W_2(\nu_1))=0.
$$

The complement of the Zariski open subset $\MO_{g-1}$ of $\M_{g-1}$
has complex codimension $g-3$ (cf. Section \ref{sec:moduli}). It
follows that $L_2$ also has complex codimension $g-3$ in $\M_{g-1,1}$.
Thus
$$
H^2(L_2; \pitil_2^\ast \W_2(\nu_1)) \cong 
H^2(\M_{g-1,1}; \W_2(\nu_1))
$$ 
when $g-3\ge 3$. The mapping class group of $\M_{g-1,1}$ is
$\Ga_{g-1,1}$, and their rational cohomology are the same. In
particular,
$$
H^2(\M_{g-1,1}; \W_2(\nu_1)) \cong H^2(\Ga_{g-1,1}; W_2(\nu_1)) =
H^2(\Ga_{g-1,1}; H^1(S;\Q))
$$
where $S$ is a reference surface of genus $g-1$. In Lemma \ref{lem:h3}
below (based on a result of Harer) we show that $H^2(\Ga_{g,1};
H^1(S;\Q))$ when $g\ge 5$. This implies that $C_\infty^{0,2}=0$, and
therefore both homomorphisms $\psitil$ and $\psi$ are zero
homomorphisms when $g\ge 6$.
\end{proof}

\begin{lemma}
\label{lem:h3}
If $g\ge 5$, then $H^2(\Ga_{g,1}; H^1(S;\Q))=0$. 
\end{lemma}

\begin{proof}
All cohomology groups are considered with rational coefficients. The
homomorphism $\Ga_{g,1}^1 \to \Ga_{g,1}$ is defined by forgetting a
fixed point, and therefore is surjective. We can choose a fixed point
in a neighborhood of the base point of a fixed tangent vector. This
determines a splitting of the homomorphism above. As the associated
spectral sequence has two rows, the existence of splitting implies
that the spectral sequence degenerates at $E_2$. Hence,
$H^2(\Ga_{g,1}; H^1(S))$ is a direct summand of $H^3(\Ga_{g,1}^1)$
\cite[Sec.~7]{harer:third}. Thus it suffices to prove that
$H^3(\Ga_{g,1}^1)=0$.

There is a short exact sequence of groups
$$
1 \to  \Z \to  \Ga_{g,2} \to  \Ga_{g,1}^1 \to  1.
$$
It determines a Gysin long exact sequence
$$
\dots \to  H^1(\Ga_{g,1}^1) \to  H^3(\Ga_{g,1}^1) \to  
H^3(\Ga_{g,2}) \to  \dots 
$$
We know that the last term is trivial according to Theorem 3.1 from
\cite{harer:fourth}. The first term is trivial by
\cite[Prop.~5.2]{hain:torelli}. It follows that the middle term
$H^3(\Ga_{g,1}^1)$ is also zero.
\end{proof}

\begin{remark}
In Theorem 3.1 from \cite{harer:fourth} Harer gives an explicit
description of a basis of $H_3(\Ga_{4,2}) \cong \Q$. Using this one
can deduce that $H_3(\Ga_{4,1}^1)$ is trivial, and therefore that
$H^2(\Ga_{4,1}; H^1(S;\Q))$ is trivial.
\end{remark}

Recall that each irreducible symplectic local system over $\M_g$ is
determined by its highest weight $\la$. If $\la_1, \dots, \la_g$ is a
set of fundamental weights of $\Symp_{2g}$, then $\la$ is uniquely
expressed as $\sum_{i=1}^g a_i \la_i$ for some non-negative integers
$a_i$. We defined $|\la|$ to be $\sum_{i=1}^g i a_i$. 

\begin{definition}
We say that an irreducible symplectic local system over $\M_g$
determined by the highest weight $\la$ is {\it even} if $|\la|$ is
even, and it is {\it odd} if $|\la|$ is odd. 
\end{definition}

The following corollary is a consequence of the proof of the main
theorem.

\begin{corollary}
\label{cor:even:main}
If $\V$ is an even local system, then the natural map
$$
IH^2(\Mtil_g;\V) \to  H^2(\M_g;\V)
$$ 
is an isomorphism when $g\ge 4$.
\end{corollary}

\begin{proof}
The estimate $g\ge 6$, rather than $g\ge 4$, appears in the proof of
Lemma \ref{lem:g6}. This lemma deals with the case when a symplectic
local system $\V$ restricted to $N^\ast$ has a direct summand
isomorphic to the irreducible local system $\pi^\ast (\W_1(0) \btimes
\W_2(\nu_1))$. Note that $\V$ contains such direct summand if and only
if the corresponding algebraic representation $V$ of $\Symp_{2g}$
restricted to the subgroup $\Spl_2 \times \Symp_{2g-2}$ contains a
copy of $W_1(0) \btimes W_2(\nu_1)$. The branching rule of
$\Symp_{2g}$ over $\Spl_2 \times \Symp_{2g-2}$ respects even and odd
components. Therefore if $V$ is even, then its restriction cannot
contain $W_1(0) \btimes W_2(\nu_1)$. This implies that in this case
$\psitil$ is trivial for all $g\ge 4$. \end{proof}


\section{Mixed Hodge theory}
\label{sec:mhs}

In this section we consider the mixed Hodge structure on $H^2(\M_g;
\V)$ where $\V$ is an irreducible symplectic local system. We prove
that the mixed Hodge structure on $H^2(\M_g;\V)$ is pure when $g\ge
6$. We also prove that if $g= 3,4,5$, then the mixed Hodge structure
on $H^2(\M_g;\V)$ has at most two weights. In this section we assume
that $g\ge 3$.

We use results of the theory of mixed Hodge modules developed by
M.~Saito. For definitions and results we refer the reader to
\cite{saito:intro}, \cite{saito:mhm}. In this paper we use only the
formal properties of mixed Hodge modules. 

\begin{notation}
Let $H=(H_\Q,H_\C,W_\bul,F^\bul)$ be a rational mixed Hodge structure
where $\W_\bul$ denotes the weight filtration, and $F^\bul$ denotes
the Hodge filtration. Denote the graded quotient $W_kH_\Q
/W_{k-1}H_\Q$ by $\Gr{W}{k}H$. We shall say that an integer $m$ is a
{\it weight} of a mixed Hodge structure $H$ if $\Gr{W}{m}H\ne 0$. We
use abbreviations: MHS for mixed Hodge structure, and MHM for mixed
Hodge module.
\end{notation}

In \cite{deligne:hodge} Deligne proved that the rational cohomology of
every quasi-projective variety possesses a natural MHS. In
\cite{saito:mhm} Saito proved that the cohomology and intersection
cohomology of an algebraic variety with coefficients in an admissible
variation of MHS carry MHSs. The definition of an admissible
variation of MHS is given for curves in \cite{steenbrink:zucker}, and
in general in \cite{kashiwara} (also see \cite[2.1]{saito:intro}).
There is a strong belief that when both MHSs of Deligne and Saito
exist they are the same.

Let $\V$ be an irreducible symplectic local system over $\M_g$
determined by highest weight $\la$. This is clear that the restriction
of the local system $\V$ to $\MO_g$ underlies a polarized variation of
Hodge structure of geometric origin. Therefore the restriction of $\V$
to $\MO_g$ is an admissible variation of Hodge structure. The local
system $\V$ is irreducible, therefore the corresponding variation of
Hodge structure is unique up to Tate twist
\cite[Prop.~8.1]{hain:torelli}. We fix $\V$ as a variation of Hodge
structure by decreeing its weight to be $|\lambda|$.

According to the theory of MHMs both $IH^q(\Mtil_g;\V)$ and $H^q (\M_g;
\V)$ carry natural MHSs \cite[pp.~146-147]{saito:intro}. The MHS on
$H^q (\M_g; \V)$ can be defined using either the smooth covers $\M_g
[l]$ for $l\ge 3$, or the isomorphism $H^q (\M_g; \V) \cong IH^q
(\M_g; \V)$ where in the second term we consider the restriction of
$\V$ to $\MO_g$. This is easy to check that all these ways lead to the
same MHS.

\begin{theorem}
\label{thm:purity}
If $g\ge 6$, or if $g\ge 4$ and $\V$ is an even local system, then the
mixed Hodge structure on $H^2(\M_g;\V(\lambda))$ is pure of weight
$2+|\lambda|$.
\end{theorem}

\begin{proof}
The theory of MHM implies that the restriction 
\begin{equation*}
IH^2(\Mtil_g;\V) \to H^2(\M_g;\V)
\end{equation*}
is a morphism of MHSs, and according to Theorem \ref{thm:main} and
Corollary \ref{cor:even:main} this is an isomorphism. The space
$\Mtil_g$ is a projective variety. It follows that the MHS on
$IH^2(\Mtil_g;\V)$ is pure of weight $2+|\lambda|$
\cite[pp.~221-222]{saito:mhm}. Thus the MHS on $H^2(\M_g;\V)$ is also
pure of weight $2+|\lambda|$.
\end{proof}

In the rest of this section we deal with the MHS on $H^2(\M_g [l];\V)$
where $\M_g [l]$ is the moduli space of curves with a level $l$
structure, and $\V$ is a symplectic local system $\V(\la)$ which
underlies a variation of Hodge structure of weight $|\la|$. We assume
that $l\ge 3$, and therefore $\M_g [l]$ is smooth and $\V$ is a
genuine (not only orbifold) local system. There exists a natural MHS
on $H^q(\M_g [l]; \V)$ for each $q\ge 0$.

\begin{theorem}
\label{thm:semipurity}
If $l\ge 3$ and $g \ge 3$, then $\Gr{W}{k}H^2(\M_g[l];\V)=0$ for $k
> 3+|\lambda|$ and $k<2+|\lambda|$.
\end{theorem}

\begin{proof}
In the beginning we recall some facts from Section \ref{sec:comp}. The
moduli space $\M_g [l]$ has the Satake compactification $\Mtil_g [l]$
which is a projective variety. The boundary $\Mtil_g [l] - \M_g [l]$
has codimension two in $\Mtil_g [l]$, and each codimension two
irreducible component $\Phi_1^\beta$ has a Zariski open subset
$Z_\beta$ such that the subsets $Z_\beta$ do not intersect each other,
and each of them is isomorphic to a smooth Zariski open subset of
$\M_1 [l] \times \M_{g-1} [l]$.

\begin{notation}
We denote by $\calS$ the intersection cohomology sheaf $\calIC(\V)$ on
$\Mtil_g [l]$. The following diagrams defines the notation for the
inclusions:
$$
\M_g [l] \stackrel{i}{\hookrightarrow} 
\M_g [l] \cup (\cup_\beta Z_\beta)
\stackrel{j}{\hookleftarrow} \cup_\beta Z_\beta,
$$
and we denote by $j_\beta$ the restriction of $j$ to $Z_\beta$. This
notation is similar to that in Section \ref{sec:main}.
\end{notation}

It follows that one has an exact sequence
$$
0 \to IH^2(\Mtil_g[l];\V) \to H^2(\M_g[l];\V) \to
H^3(\cup_\beta Z_{\beta}; j^! \calS)
$$
in the category of MHSs. Taking graded quotients with respect to
weight filtration is an exact functor. Therefore for every $k$ there
is an exact sequence
$$
0 \to \Gr{W}{k}IH^2(\Mtil_g[l];\V) \to \Gr{W}{k}H^2(\M_g[l];\V) 
\to \Gr{W}{k}H^3(\cup_\beta Z_{\beta}; j^! \calS).
$$ 
Since the space $\Mtil_g[l]$ is a projective variety, and $\V$ is a
polarized variation of Hodge structure of geometric origin of weight
$|\la|$, the intersection cohomology $IH^2(\Mtil_g[l];\V)$ has a pure
MHS of weight $2+|\la|$. To prove the theorem we will show that
$\Gr{W}{k}H^3(\cup_\beta Z_{\beta}; j^! \calS)= 0$ unless $k=3+
|\la|$.

As the sets $Z_\beta$ are disjoint it suffices to show that each
$H^3(Z_\beta; j_\beta^! \calS)$ has a pure MHS of weight
$3+|\la|$. From now on we fix an arbitrary index $\beta$, and omit
$\beta$ from the notation for $Z_\beta$ and $j_\beta$.

The sheaf $j^! \calS$ is constructible, and $\Mtil_g [l]$ is
equi-singular along $Z$.  Therefore $\calH^3 j^! \calS$
is a local system over $Z$. The standard argument implies that there
is an isomorphism of MHSs
\begin{equation*}
H^3(Z; j^! \calS) \cong H^0 (Z; \calH^3 j^! \calS),
\end{equation*}
and there is an isomorphism of MHMs
\begin{equation}
\label{sheaves}
\calH^3 j^! \calS \cong \calH^2 j^\ast i_\ast \V.
\end{equation}
We will show that these MHMs are pure of weight $3+|\la|$.

Recall that $j^\ast i_\ast \V$ expresses cohomology of the link of $Z$
in $\M_g [l] \cup Z$. The inverse image of $Z$ under the birational
morphism $\alpha^l: \Mbar_g [l] \to \Mtil_g [l]$ is a smooth locally
closed divisor. We denote it by $Y$. Then the link of $Z$ in $\M_g [l]
\cup Z$ is the same as the link of $Y$ in $\M_g [l] \cup Y$. We use
this to find the weights on $\calH^2 j^\ast i_\ast \V$.

The following commutative diagram introduces the notation: 
$$
\begin{CD}
\M_g[l]     @>\kappa >> \M_g[l] \cup Y       @<\mu << Y  \\
@V=VV          @V\alpha^l VV                @V\pibar VV  \\
\M_g[l]   @>i >>  \M_g[l] \cup Z @<j <<    Z.   \\  
\end{CD}
$$
The local link cohomology functor of $Y$ is $\mu^\ast \kappa_\ast\V$.
Therefore one expects that $j^\ast i_{\ast}\V \simeq \pibar_\ast
\mu^\ast \kappa_\ast\V$. (The sign $\simeq$ denotes an isomorphism in
the derived category of MHMs.) Indeed, both $\alpha^l$ and $\pibar$ are
proper maps, therefore $\alpha^l_\ast= \alpha^l_!$ and $\pibar_\ast=
\pibar_!$. It follows that for an arbitrary sheaf $\calF$ on $\M_g [l]
\cup Y_\beta$ one has that $j^\ast\alpha^l_\ast \calF \simeq \pibar_\ast
\mu^\ast \calF$ \cite[Prop.~10.7]{borel:ic}. Therefore
\begin{equation}
\label{links}
j^\ast i_{\ast}\V\simeq j^\ast\alpha^l_\ast \kappa_\ast\V
\simeq \pibar_\ast \mu^\ast \kappa_\ast\V.
\end{equation}
Thus $\calH^2 j^\ast i_{\ast} \V \cong \calH^2 \pibar_\ast \mu^\ast
\kappa_\ast\V$ is an isomorphism of MHMs.

The variation of Hodge structure $\V$ on $\M_g [l]$ extends to a
variation of Hodge structure on $\M_g [l] \cup Y_\beta$ because $\V$
is pulled back from $\A_g[l]$. We denote its restriction to $Y$ by
$\Vbar$. Then $\mu^\ast \kappa_\ast\V \simeq \mu^\ast \kappa_\ast\Q
\otimes \Vbar$ where $\Q$ denotes the constant variation of Hodge
structure of weight zero with the fiber isomorphic to $\Q$.

Denote by $D\calF$ the dual of $\calF$ in the derived category of
MHMs. The spaces $\M_g[l] \cup Y$ and $Y$ are smooth, therefore $D\Q
\simeq \Q[2n](n)$ and $D\Q_{Y_\beta} \simeq \Q_{Y_\beta}[2n-2](n-1)$.
It follows that there is a string of isomorphisms in the derived
category of MHMs:
%
\begin{multline*}
\mu^!\Q \simeq  D_{Y_\beta}(\mu^\ast D\Q)\simeq 
D_{Y_\beta}(\mu^\ast \Q[2n](n))  \simeq   
D_{Y_\beta}(\Q_{Y_\beta}[2n](n)) \\  \simeq 
(D_{Y_\beta}\Q_{Y_\beta})[-2n](-n) \simeq  \Q_{Y_\beta}[-2](-1).
\end{multline*}
%
Using this and the distinguished triangle
\ecouple{\mu^!\Q} {}
	{\Q} {}
	{\mu^\ast \kappa_\ast\Q.} {[1]}
one can deduce that 
$$ \calH^0 \mu^\ast \kappa_\ast\Q \cong \calH^0 \Q, \quad
\calH^1 \mu^\ast \kappa_\ast\Q \cong \calH^2 \mu^!\Q, \quad 
\text{and} \quad 
\calH^q \mu^\ast \kappa_\ast\Q = 0
$$ 
for $q\ge 2$. It follows that $\calH^0 \mu^\ast \kappa_\ast\V$ is a
pure Hodge module of weight $|\la|$, and $\calH^1 \mu^\ast
\kappa_\ast\V$ is a pure Hodge module of weight $2+|\la|$.

According to \cite[1.20]{saito:intro}, there is a (perverse) spectral
sequence in the category of MHMs:
\begin{equation*}
E_2^{p,q}= \calH^p \pibar_\ast(\calH^q\mu^\ast \kappa_\ast\V)
\Rightarrow \calH^{p+q}\pibar_\ast \mu^\ast \kappa_\ast\V.
\end{equation*}
As all spaces involved are smooth the perverse spectral sequence
coincides with the ordinary one. It has only two non-zero rows. Thus
there is an exact sequence of MHMs
$$
\calH^2 \pibar_\ast (\calH^0 \mu^\ast \kappa_\ast\V) \to
\calH^2 \pibar_\ast \mu^\ast \kappa_\ast\V \to
\calH^1 \pibar_\ast (\calH^1\mu^\ast \kappa_\ast\V).
$$
The map $\pibar$ is proper. Therefore $\calH^2\pibar_\ast (\calH^0
\mu^\ast \kappa_\ast\V)$ is pure of weight $2+|\la|$, and $\calH^1
\pibar_\ast (\calH^1 \mu^\ast \kappa_\ast\V)$ is pure of weight
$3+|\la|$. Consequently, we have that
$$
\Gr{W}{k} \calH^2 \pibar_\ast \mu^\ast \kappa_\ast\V=0
$$ 
for $k> 3+|\la|$ and $k< 2+|\la|$.

Since $\V$ is a variation of Hodge structure of geometric origin of
weight $|\la|$, the intersection cohomology sheaf $\calS$ underlies a
pure Hodge module of weight $|\la|$. Therefore $j^!\calS$ is a MHM of
weight $\ge |\la|$ \cite[Prop.~1.7]{saito:intro}. It follows that
$\Gr{W}{k} \calH^3 j^!\calS=0$ for $k< 3+|\la|$. Combining the last
two paragraphs, and isomorphisms \eqref{sheaves}, \eqref{links} one
gets that $\calH^3 j^! \calS$ is pure of weight $3+|\la|$. 
\end{proof}

\begin{corollary}
\label{cor:semipurity}
Let $\V(\la)$ be a symplectic local system over $\M_g$ underlying a
variation of Hodge structure of weight $|\la|$. If $g\ge 3$, then
$\Gr{W}{k} H^2(\M_g; \V)=0$ for $k> 3+|\la|$ and $k<2+|\la|$.
\end{corollary}

\begin{proof}
Choose $l\ge 3$. One has as isomorphism 
$$
H^2(\M_g; \V) \cong H^2(\M_g[l]; \V)^{\Symp_{2g}(\Z/l\Z)} 
$$
in the category of MHSs. The weights of the right hand side are
$2+|\la|$ and $3+|\la|$ according to the theorem above. Therefore the
same is true for the left hand side. 
\end{proof}



\begin{thebibliography}{99}

\bibitem{acgh:II}
E.~Arbarello, M.~Cornalba, F.~Griffiths, J.~Harris,
{\it Geometry of Algebraic Curves II},
preliminary manuscript, 1994.

\bibitem{baily:1962}
W.~Baily,
On the theory of theta-functions, the moduli of 
abelian varieties and the moduli of curves, 
{\it Ann. Math.} (2) 75 (1962), 342--381.

\bibitem{bbd}
A.~Beilinson, J.~Bernstein, P.~Deligne,
Faisceaux Pervers, 
{\it Ast\'erisque} 100, 1983. 

\bibitem{birman:braids}
J.~Birman,
{\it Braids, Links and Mapping Class Groups}, 
Ann. of Math. Studies 82, Princeton University Press, 
Princeton, 1982. 

\bibitem{borel:ic}
A.~Borel et al.,
{\it Intersection Cohomology}, Progress in 
Math. 50, Birkh\"{a}user, Boston, 1984.

\bibitem{bredon:sheaf}
G.~Bredon,
{\it Sheaf Theory}, 
McGraw--Hill Book Company, 1967.

\bibitem{deligne:hodge}
P.~Deligne,
Th\'eorie de Hodge I, {\it Actes, Congr\`es Intern. Math.}, 
Nice, 1970, 425--430. 
Th\'eorie de Hodge II, {\it Publ. Math. I.H.E.S.} 40 (1971), 5--58. 
Th\'eorie de Hodge III, {\it Publ. Math. I.H.E.S.} 44
(1974), 5--77.

\bibitem{deligne:mumford}
P.~Deligne, D.~Mumford, 
The irreducibility of the space of curves of given genus,
{\it Publ. Math. Inst. Hautes \'Etudes Sci.} 36 (1969), 75--110.

\bibitem{durfee:saito}
A.~Durfee, M.~Saito,
Mixed Hodge structures on the intersection cohomology of links,
{\it Compositio Math.} 76 (1990), 49--67.  

\bibitem{eisenbud:harris}
D.~Eisenbud, J.~Harris,
Progress in the theory of complex algebraic curves, 
{\it Bulletin (New Series) of AMS} 21, 2, Oct. (1989), 205--232. 

\bibitem{fulton:harris}
W.~Fulton, J.~Harris,
{\it Representation Theory, a First Course}, 
GTM 129, Springer--Verlag, 1991.

\bibitem{gmII}
M.~Goresky, R.~MacPherson,
Intersection homology II, {\it Invent. Math.} 
71 (1983), 77--129.

\bibitem{hain:completion}
R.~Hain,
Completions of mapping class groups and the cycle $C-C^-$, 
{\it Mapping class groups and moduli spaces of Riemann surfaces} 
(C.~B\" odigheimer and R.~Hain, eds.), Contemp. Math. 
150, AMS, 1993, 75--105.

\bibitem{hain:torelli}
R.~Hain,
Torelli groups and geometry of moduli spaces of curves, 
{\it Current Topics in Complex Algebraic Geometry} 
(C.~H.~Clemens and J.~Kollar, eds.), MSRI publications 28, 
Cambridge University Press, 1995, 97--143. 

\bibitem{hain:lietor}
R.~Hain, 
Infinitesimal presentations of the Torelli group, 
preprint, 1995. 
(Available from {\tt http://www.math.duke.edu/faculty/hain/})

\bibitem{harer:second}
J.~Harer,
The second homology group of the mapping class group of 
an orientable surface, 
{\it Invent. Math.} 72 (1983), 221--239. 

\bibitem{harer:stability}
J.~Harer, 
Stability of the homology of the mapping class groups 
of orientable surfaces, 
{\it Ann. of Math.} (2) 121 (1985), 215--249. 

\bibitem{harer:review}
J.~Harer,
The cohomology of the moduli space of curves, 
{\it CIME notes, Theory of moduli} (E.~Sernesi, ed.),
Lecture Notes in Math. 1337,
Springer--Verlag, Berlin and New York, 1988, 138--221.

\bibitem{harer:third}
J.~Harer,
The third homology group of the moduli space of curves,
{\it Duke Math. J.} 63 (1991), 25--55. 

\bibitem{harer:fourth}
J.~Harer, 
The fourth homology group of the moduli space of curves,
preprint, 1993. 

\bibitem{ivanov:stability}
N.~Ivanov,
On stabilization of the homology of Teichm\"uller 
modular groups,
{\it Algebra i Analyz} 1 (1989), 110-126;
English Translation: {\it Leningrad J. of Math.} 1 
(1990), 675--691. 

\bibitem{ivanov:twist}
N.~Ivanov,
On the homology stability for Teichm\" uller modular groups. 
Closed surfaces and twisted coefficients, 
{\it Mapping class groups and moduli spaces of Riemann surfaces}
(C.~B\" odigheimer and R.~Hain, eds.), Contemp. Math. 150, 
AMS, 1993, 149--194. 

\bibitem{johnson}
D.~Johnson,
The structure of the Torelli group I: a finite set of 
generators for $\mathcal I$, 
{\it Ann. of Math.} 118 (1983), 423--442. 

\bibitem{kashiwara}
M.~Kashiwara,
A study of variation of mixed Hodge structure, 
{\it Publ. RIMS} 22 (1986), 991-1024. 

\bibitem{knudsen:II}
F.~Knudsen,
The projectivity of the moduli space of stable curves, II:
the stacks $M_{g,n}$,
{\it Math. Scandinavica} 52 (1983), 161--199.

\bibitem{knudsen:III}
F.~Knudsen,
The projectivity of the moduli space of stable curves, III:
the line bundles on $M_{g,n}$, and a proof of the projectivity 
of $\bar M_{g,n}$ in characteristic 0, 
{\it Math. Scandinavica} 52 (1983), 200--212. 

\bibitem{looijenga:prym}
E.~Looijenga, 
Smooth Deligne--Mumford compactifications by means of Prym level
structures,
{\it J. Algebraic Geometry} 3 (1994), 283--293.

\bibitem{looijenga:stable}
E.~Looijenga,
Stable cohomology of the mapping class group with 
symplectic coefficients and of the Abel--Jacobi map, 
{\it J. Algebraic Geometry} 5 (1996), 135--150.

\bibitem{mostafa}
S.~Mostafa,
Die Singularit\"aten der Modulmannigfaltigkeit $\Mbar_g (n)$
der stabilen Kurven vom Geschlecht $g\ge 2$ mit 
$n$-Teilungspunktstruktur,
{\it J. Reine Angew. Math.} 343 (1983), 81--98.

\bibitem{mumford:git}
D.~Mumford,
{\it Geometric Invariant Theory}, 
Springer--Verlag, Berlin Heidelberg, 1965. 

\bibitem{mumford:stable}
D.~Mumford,
Stability of projective varieties,
{\it L'Enseignement Math.} II Ser., XXIII-fasc. 1--2 (1977), 39--110.

\bibitem{mumford:red}
D.~Mumford, 
{\it The Red Book of Varieties and Schemes},
Springer--Verlag, Berlin Heidelberg New York, 1988. 

\bibitem{oort:steenbrink}
F.~Oort, J.~Steenbrink,
The local Torelli problem for algebraic curves, 
{\it Journ\'ees de g\'eom\'etrie alg\'ebrique d'Angers} 
(A.~Beauville, ed.), 1979, 157--204. 

\bibitem{pikaart}
M.~Pikaart,
An orbifold partition of $\overline{M}_g^n$, 
{\it The moduli space of curves} 
(R.~Dijkgraaf, C.~Faber, G.~van der Geer, eds.), 
Birkh\"auser, 1995, 467--482.

\bibitem{jong:pikaart}
M.~Pikaart, A.~de Jong,
Moduli of curves with non-abelian level structure,
{\it The moduli space of curves} 
(R.~Dijkgraaf, C.~Faber, G.~van der Geer, eds.), 
Birkh\"auser, 1995, 483--509. 

\bibitem{popp}
H.~Popp,
{\it Moduli Theory and Classification Theory of Algebraic Varieties},
Lecture Notes in Math. 620, Springer--Verlag, Berlin Heidelberg,
1977.

\bibitem{saito:intro}
M.~Saito,
Introduction to mixed Hodge modules, 
{\it Th\'eorie de Hodge}, Luminy, Juin 1987, 
{\it Ast\'erisque} 179-180, 1989, 145-162.

\bibitem{saito:mhm}
M.~Saito,
Mixed Hodge modules, 
{\it Publ. RIMS} 26 (1990), 221--335. 

\bibitem{satake}
I.~Satake,
On the compactification of the Siegel space,
{\it Indian Math. Soc.} 20 (1956), 259--281.

\bibitem{steenbrink:zucker}
J.~Steenbrink, S.~Zucker, 
Variation of mixed Hodge structure I, 
{\it Invent. Math.} 80 (1985), 489-542. 

\end{thebibliography}
\end{document}